\documentclass[reqno, 12pt]{amsart}

\def\C{{\mathbb C}}

\def\Q{{\mathbb Q}}
\def\Qp{{\mathbb Q}_p}
\def\Z{{\mathbb Z}}
\def\Zp{{\mathbb Z}_p}
\def\R{{\mathbb R}}
\def\A{{\mathbf A}}
\def\K{{\mathbf K}}
\def\k{{\mathbf k}}
\def\bfC{{\mathbf C}}
\def\KV{\K_{\mathbf V}}

\def\Fp{{{\mathbb F}_p}}
\def\F{{\mathbb F}}
\def\Fq{{{\mathbb F}_q}}
\def\Fr{{{\mathbb F}_r}}

\newtheorem{cor}{Corollary}
\newtheorem{prop}{Proposition}

\theoremstyle{definition}
   
\newtheorem{conj}{Conjecture}

\theoremstyle{remark}
\newtheorem{rem}{Remark}        
\newtheorem{rems}{Remarks}      
\newtheorem{example}{Example}

\begin{document}

\title[A Riemann Hypothesis for characteristic $p$ $L$-functions]
{A Riemann Hypothesis for characteristic $p$ $L$-functions}
\author{David Goss}
\thanks{This paper is respectfully dedicated to {\sc Bernard Alter}
and {\sc Shirley Hasnas}}
\address{Department of Mathematics\\ The Ohio State University\\ 231 W.
$18^{\text{th}}$ Ave. \\ Columbus, Ohio 43210}
\email{goss@math.ohio-state.edu}
\date{October, 1999}

\begin{abstract}
We propose analogs of the classical Generalized Riemann Hypothesis and
the
Generalized Simplicity Conjecture for the characteristic $p$ $L$-series 
associated
to function fields over a finite field. These analogs are based on the
use of absolute values. Further we use absolute values to give similar
reformulations of the classical conjectures (with, perhaps, finitely many
exceptional zeroes). We show how both sets of
conjectures behave in remarkably similar ways.
\end{abstract}

\maketitle


\section{Introduction}\label{intro}
The arithmetic of function fields attempts to create a model of classical
arithmetic using Drinfeld modules and related constructions such
as shtuka, $\A$-modules, $\tau$-sheaves, etc. Let $\k$ be one such
function field over a finite field $\Fr$ and let $\infty$ be a fixed
place of $\k$ with completion $\K=\k_\infty$. It is well known
that the algebraic closure of $\K$ is infinite dimensional over $\K$
and that, moreover, $\K$ may have infinitely many distinct extensions of
a bounded degree. Thus function fields are inherently ``looser'' than
number fields where the fact that $[\C\colon \R]=2$ offers
considerable restraint. As such, objects of classical number theory
may have many different function field analogs.

Classifying the different aspects of function field arithmetic
is a lengthy job. One finds for instance that there are two
distinct analogs of classical $L$-series. One analog comes from
the $L$-series of Drinfeld modules etc., and is the one of interest
here. The other analog arises from the $L$-series of modular forms
on the Drinfeld rigid spaces, (see, for instance, \cite{go2}).
It is a very curious phenomenon that the first analog possesses no
obvious functional equation whereas the second one indeed has
a functional equation very similar to the classical versions.
It is even more curious that the $L$-series of Drinfeld
modules and the like seem to possess the correct analogs of the {\it
Generalized Riemann Hypothesis} and the {\it Generalized Simplicity Conjecture}
(see Conjecture \ref{clsimpl} below). It is our purpose here to 
define these characteristic $p$ conjectures and show just how close they
are to their classical brethren.

That there might be a good Riemann Hypothesis in the characteristic $p$
theory first arose from the ground-breaking work \cite{w1} of
Daqing Wan. In this paper, and in the simplest possible case, Wan
computed the valuations of 
zeroes of an analog of the Riemann zeta function via the technique
of Newton polygons. This immediately implied
that these zeroes are all simple and lie on a ``line.'' However,
because of the great size of the function field arena (as mentioned
above), it was
not immediately clear how to then go on to
state a Riemann Hypothesis in the function field
case which worked for all places of $\k$ (as explained in this paper)
and all functions arising from arithmetic.

Recently, the $L$-functions of function field arithmetic were
analytically continued in total generality (as general as one could
imagine from the analogy with classical motives). This is due to the
forthcoming work of G.\ Boeckle and R.\ Pink 
\cite{bp1} where an appropriate cohomology
theory is created. This theory, combined with certain estimates provided
by Boeckle, on the one hand, and Y.\ Amice \cite{am1}, on the other, actually
allows
one to analytically continue the non-Archimedean measures
associated to the $L$-series; the analytic continuation of the
$L$-series themselves then arises as a corollary. In particular we deduce
that {\it all} such $L$-functions, viewed at all places of $\k$,
have remarkably similar analytic properties (for instance, their expansion
coefficients all decay exponentially --- see the discussion after
Remarks \ref{notapply}).

Motivated by these results, we re-examined the work of Wan
and those who came after him (\cite{dv1}, \cite{sh1}). In seeking to
rephrase Wan's results in such a way as to avoid having to compute
Newton polygons (which looks to be exceedingly complicated in general), we
arrived at a statement involving only the use of {\it absolute values} of
zeroes (as opposed to the absolute values of expansion coefficients which are
used in Newton polygons). 
The use of absolute values in phrasing such a possible Riemann
Hypothesis seems to be very fruitful. For instance, it offers
a unification with local Riemann Hypotheses (which are always formulated
in terms of absolute values of the zeroes). More strikingly, it also
suggests a suitable reformulation of the classical GRH (with, perhaps,
finitely many exceptional zeroes) as well
as the simplicity conjectures (see Conjecture \ref{classGRH3} and Proposition
\ref{c5.1}). Finally, as explained after Remarks \ref{feuse}, the conjectures
presented here go a very long way towards explaining the lack of a classical
style functional equation associated to the $L$-series of Drinfeld modules
etc.

Upon examining these new ``absolute value conjectures'' 
in both theories, one finds
that they behave remarkably alike. So much so that they seem to
almost be two instances of one Platonic mold. This certainly adds to our
sense that the function field statements may indeed be the correct
ones. Moreover, because the algebraic closure of $\K$ is so vast and
contains inseparable extensions, the
function field theory offers insight into these statements not available
in number fields. For instance, due to the existence of 
inseparable extensions, one needs both
the function field analog of the GRH {\it and} the function field
analog of the Simplicity Conjecture (Conjecture \ref{charpsim})
to truly deduce that the zeroes
(or almost all of them) lie on a line! Because $\C$ is obviously
separable over $\R$ one only needs the GRH, (reformulated as
Conjecture \ref{classGRH3}) classically.

It should be noted that we do not yet know the implications of our
function field conjectures. However, it is our hope that such information
will be found as a byproduct of the search for a proof of them.
Moreover, because of the strong formal analogies between the number field
and function field conjectures as 
presented here, we believe that insight obtained 
for one type of global field may lead to insight for the other.


\section{Basic Statements}\label{basic}
In this section we present the statements of the characteristic $p$ Generalized
Riemann Hypothesis as well as the characteristic $p$  Generalized Simplicity
Conjecture. 
We will begin by recalling the classical versions of
these conjectures. We will work with
classical abelian characters over a number field
$\mathfrak L$. However, the
reader will easily see the simple modifications necessary to handle other
classical $L$-series. Moreover, of course, if $\chi$ is one such
abelian character,
then the analytic continuation of $L(\chi,s)$ has long been known. Let
$\Lambda(\chi,s)$ be the completed $L$-series (so $\Lambda(\chi,s)$ contains
the $\Gamma$-factors at the infinite primes). One knows that there is a
functional equation relating
 $\Lambda (\chi,s)$ and $\Lambda (\bar{\chi},1-s)$.

Following Riemann's original paper, we define
$$\Xi(\chi,t):=\Lambda (\chi,1/2+it)\,.$$
It is very easy to see that the functional equation for $\Lambda$ translates
into one for $\Xi$ relating $\Xi(\chi,t)$ and
$\Xi(\bar{\chi},-t)$.
\smallskip

\begin{conj}\label{classGRH1}
The zeroes of $\Lambda (\chi,s)$ lie on the line $\{s=1/2+it\mid t\in \R\}$.
\end{conj}
\smallskip

\noindent
Conjecture \ref{classGRH1} is obviously the {\it Generalized Riemann
Hypothesis} for abelian $L$-series. It
may clearly be reformulated (as Riemann did)
in the following way.
\smallskip

\begin{conj}\label{classGRH2}
The zeroes of $\Xi(\chi,t)$ are real.
\end{conj}
\smallskip

All known zeroes of the Riemann zeta function, $\zeta(s)=L(\chi_0,s)$ where
$\chi_0$ is the trivial character and $\mathfrak{L}=\Q$, 
have been found to be simple. This is
codified in the following conjecture which the author learned from
J.-P. Serre. In it we {\it explicitly} assume that our number field 
$\mathfrak L$ is now the base field $\Q$.\smallskip

\begin{conj}\label{clsimpl}
\noindent 1.  The zeroes of $\Lambda(\chi ,s)$ 
should be simple.

\noindent 2.  $s=1/2$ should be a zero of $\Lambda(\chi ,s)$ only 
if $\chi$ is real and the functional equation of $\Lambda (\chi ,s)$ has 
a minus sign.

\noindent 3.  If $\chi$ and $\chi '$ are distinct, then the 
zeroes of $\Lambda (\chi ,s)$ not equal to $s=1/2$ should be distinct 
from those of $\Lambda(\chi ',s)$.
\end{conj}
\smallskip

\noindent
We shall call Conjecture \ref{clsimpl} the {\it Generalized Simplicity
Conjecture} (``GSC''). It is in fact expected to hold for general (not necessarily
abelian) $L$-series
(see e.g., Conjecture 8.24.1 of \cite{go1}).

We turn next to the function field versions of the above
conjectures.
As this theory is certainly not as well known as its
classical counterparts, we begin by reviewing it. For more the reader
can consult \cite{go1}; for a short and very readable introduction to
Drinfeld modules the reader may consult \cite{h1}. 
Let $\mathcal X$ be a smooth projective
geometrically connected curve over the finite field $\Fr$, $r=p^m$, and
let $\infty\in \mathcal{X}$ be a fixed closed point. Let $\A$ be the
affine ring of $\mathcal{X}\backslash \infty$; so $\A$ is a Dedekind
domain with finite class group and unit group equal to $\Fr^\ast$. We
let $\k$ be the function field of $\mathcal{X}$ (= the quotient field
of $\A$) and we let $\K=\k_\infty$. Finally we let $\bfC_\infty$ be the
completion of a fixed algebraic closure of $\K$ equipped with its
canonical topology.

A particular instance is $\mathcal{X}=\mathbb{P}^1_\Fr$ and $\A=\Fr[T]$, etc.
The reader may wish to first read this paper with only this instance in
mind as the jump to general $\A$ is largely technical.

The ``standard analogy'' is that $\A$ plays the role classically
played by $\Z$, $\k$ the role of $\Q$, $\K=\k_\infty$ the role
of $\R$ and $\bfC_\infty$ the role of $\C$. We will have more to say
on this below (see Remark \ref{twoT}).

As mentioned above, the basic ``arithmetic'' objects in the
characteristic $p$ theory are Drinfeld
modules (originally presented in \cite{dr1}) and their various generalizations such as $\A$-modules. These
arise in the following manner. Let
$L$ be some field over $\Fr$ and let $\mathbb{G}^n_a$ be the $n$-th
Cartesian product of $\mathbb{G}_a$
over $L$. Let $x=\left( \begin{array}{c} x_1\\ \vdots\\x_n
\end{array}\right)\in \mathbb{G}^n_a$ and set
$$\tau (x):=\left( \begin{array}{c} x_1^r\\ \vdots\\x_n^r
\end{array}\right)\,;$$
one forms $\tau^i$ by composition for all non-negative integers
$i$.
Let $P\colon \mathbb{G}_a^n\to \mathbb{G}_a^n$ be a morphism of algebraic
groups over $L$
that is $\Fr$-linear on the geometric points. It is elementary to see
that there is a ``polynomial'' $P(\tau):=\sum_{i=0}^m a_i \tau^i$, where the
$a_i$ are $n\times n$ matrices with coefficients in $L$, such that
$$P(x)=\sum_{i=0}^m a_i \tau^i(x)=P(\tau)(x)\,.$$ 
We denote the $\Fr$-algebra (under composition) of such maps by
${\rm End}_{L}^r (\mathbb{G}_a^n)$. The mapping which takes
$P(\tau)$ to $a_0$ is readily checked to be a map of $\Fr$-algebras; we
denote $a_0$ by $P^\prime$ or  $P^\prime (\tau)$.

Suppose now that $L$ is also 
an ``$\A$-field;'' that is, there is an $\Fr$-algebra
map $\imath \colon \A\to L$. Following Anderson \cite{a1} an
{\it $\A$-module}, $\phi$, is then an $\Fr$-algebra map
from $\A$ to ${\rm End}_{L}^r(\mathbb{G}_a^n)$, $a\in \A \mapsto \phi_a$,
subject to the condition that $N_a:=\phi_a^\prime -\imath(a)\cdot I_n$ is
nilpotent. The {\it dimension} of $\phi$ is $n$.
Notice that if $n=1$, then $N_a=0$ for all $a\in \A$.
A {\it Drinfeld module} is a $1$-dimensional $\A$-module such that, for
some $a\in \A$, $\phi_a(\tau)-\imath(a)\tau^0$ is non-trivial. (In other
words, the ``trivial'' action $\psi_a(\tau):=\imath(a)\tau^0$ is not
a Drinfeld module.)  

\begin{example} \label{carlitz}
Let $\A=\Fr[T]$, $\k=\Fr(T)$, etc. The {\it Carlitz module} $C$ 
(originally presented in 1935 in \cite{c1}) is
the Drinfeld module defined over $\k$ with
$$C_T(\tau):=T\tau^0+\tau\,.$$
$C$ is the simplest, and most basic, of all Drinfeld modules.
\end{example}

Let $M:={\rm Hom}_{L}^r(\mathbb{G}_a^n,\mathbb{G}_a)$ the $\Fr$-linear
morphisms of affine algebraic groups defined over $L$.
The {\it $\A$-motive} \cite{a1} associated to an $\A$-module
$\phi$ is $M_\phi=M$
viewed as an $L\otimes_\Fr\A$-module as follows: 
Let $l\in L$, $a\in \A$, and
$m(x)\in M$. Then we set
$$(l\otimes a)\cdot m(x)=l\cdot m(\phi_a(\tau)(x))\in \mathbb{G}_a\,.$$

Passing from the module to the motive is extremely useful in the 
theory, see Section 5 of \cite{go1}. The
{\it $\tau$-sheaves} are natural generalizations of $\A$-motives.

\begin{rem} \label{twoT}
Suppose that $L$ is a field over of $\k$ and let $T\in \A$ be
non-constant. Then $T$ plays two roles in the theory (``two $T$'s''):
In $L$ the element $T$ is a scalar whereas in $\A$ one knows that
$T$ is an operator (via some
$\A$-module etc.). This is completely similar to the fact that an
integer $n$ plays two similar roles for elliptic curves over $\Q$.
The ``standard caveat'' is that it is obviously
impossible to separate the two distinct actions of an integer $n$ via a module
over ``$\Z\otimes \Z$'' as in the function field theory. 
\end{rem}

\noindent
It is important to keep the two
actions of an element $T$ separate. To do so we follow \cite{a1}
and use a notational device:
Let $A$, $k$, etc., be another copy of the basic algebras constructed
above. There is an obvious isomorphism $\theta$
from ``bold'' to ``non-bold''
making the non-bold rings $\A$-algebras. When $\A=\Fr[T]$, it is
customary to set $\theta:=\theta(T)\in A$.
The elements of the bold algebras will be the operators while
the elements of the non-bold algebras will be the scalars. In this set-up,
$\A$-modules etc., will always be defined over the non-bold scalars.

\begin{example} \label{2carlitz}
Let $C$ be the Carlitz module as in Example \ref{carlitz}.
Using the notation just introduced, the Carlitz module is the
Drinfeld module defined over $\Fr(\theta)$ with
$$C_T(\tau):=\theta\tau^0+\tau\,.$$
\end{example}

In particular, if $T\in \A$ and $x\in \mathbb{G}_a^n$, then the
action $T\cdot x$ is now unambiguously given and one may suppress the
use of $\phi$ etc.
As a bonus, the ``bold, non-bold'' notation also 
 provides an extremely useful way of classifying the
constructions of function field arithmetic.
For instance, the periods of Drinfeld modules
or $\A$-modules are scalars; thus all constructions of $\Gamma$-functions
take values in non-bold algebras,  or are ``scalar-valued.'' Similarly,
the $L$-series of modular forms are scalar-valued. On the other hand,
the $L$-series of Drinfeld modules, etc., are derived from Tate-modules
exactly as in classical arithmetic. 
As the Tate modules are, by definition, modules
over operator algebras, we see that these $L$-series are ``operator-valued.''
We now describe these operator-valued $L$-series in some detail.

In all that follows, the reader should keep in mind the 
the $L$-series of an abelian variety over a number field or
$L(\chi,s)$ mentioned above. The following
definitions are given with the goal of making the characteristic
$p$ $L$-series as close as possible to these classical $L$-series.

We begin by explaining what is meant by a ``Dirichlet series'' in the
characteristic $p$ context. 
Let $S_\infty:={\bfC}_\infty\times \Zp$; we view $S_\infty$
as a topological abelian group whose operation is written additively.
The space $S_\infty$ will supply the ``$s$'' in ``$I^s$.'' Let
$\F_\infty\subset \K$ be the constant field at $\infty$.
Let ${\rm sgn}\colon \K^\ast\to \F_\infty^\ast$ be a sign function; that is
a morphism which is the identity on $\F_\infty^\ast$. An element $z\in
\K^\ast$ is said to be {\it positive} (or {\it monic}) if and only
${\rm sgn}(z)=1$; it is very elementary to see that the 
group $\mathcal{ P}_+$ of
principal and positively generated $\A$-fractional ideals 
of $\k$ is of finite index in the total group $\mathcal{I}$ of
$\A$-fractional ideals. Let $\pi\in \K$ be a positive uniformizer. Let
$d_\infty$ be the degree of $\infty$ over $\Fr$ and set
$d(z):=-d_\infty v_\infty (z)$ where $v_\infty (\pi)=1$. If $a\in \A$ is
positive, then $d(a)$ is the degree of the finite part of the divisor
of $a$. We then put
$$\langle z \rangle=\langle z \rangle_\pi:=\pi^{-v_\infty(z)}\cdot z\,.$$
Notice that $\langle z\rangle\in U_1$ where $U_1\subset \K^\ast$ is the
subgroup of $1$-units. It is easy to see, via the binomial theorem, that
$U_1$ is a $\Zp$-module. 

Now let $I$ be any $\A$-fractional ideal and let $e$ be the
index of $\mathcal{P}_+$ in $\mathcal I$.
Thus $I^e=(i)$ for positive $i \in \A$.
Let ${\mathbf U}_1\subset \bfC_\infty$ be the group of $1$-units. Notice
that $\mathbf{U}_1$ is a divisible group (in fact,
a $\Qp$-vector space); in particular,
$\mathbf{U}_1$ is injective. Thus, if we
set 
$$\langle I\rangle := \langle i\rangle^{1/e}\,,$$
where we take the solitary root in $\mathbf{U}_1$, we obtain
a unique morphism $\langle ~\rangle\colon \mathcal{I}\to \mathbf{U}_1$ 
extending the morphism defined on $\mathcal{P}_+$ 
with $\langle (a)\rangle=\langle a \rangle$ for $a$ positive. Finally, if
$s=(x,y)\in S_\infty$ then we set
$$I^s:=x^{\deg I}\cdot \langle I\rangle^y\,.$$
If, for instance, $I=(i)$, where $i$ is positive, then
$$(i)^s=x^{d(i)}\langle i\rangle^y\,.$$ 

Since $U_1$ is a $\Zp$-module, the values $\langle I\rangle$, for $I$
an $\A$-fractional ideal, are totally inseparable over
$\K$. They generate a finite, totally inseparable, extension
denoted by $\KV$; we call $\KV$ the {\it local value field}. Of course,
$\KV=\K$ when $\A$ has class number $1$. 

A Dirichlet series $L(s)$ is then a sum
$L(s)=\sum_I c(I)I^{-s}$ 
where we sum over the ideals of $\A$ and where the
elements $c(I)$ lie in a finite extension of 
$\K$ (in practice, actually a finite extension of $\k$). Notice that if 
we set $y=0$  and $u=x^{-1}$ in $I^s$,
we obtain a characteristic $p$ version of the classical
power series arising from Artin-Weil $L$-series of function fields.
Such characteristic $p$ Dirichlet series may (and usually do) arise
from Euler products over the finite primes in an obvious sense.
Let $s=(x,y)$ be as above. By definition we have
$$L(s)=\sum_i c(I) x^{-\deg I}\langle I \rangle^{-y}\,.$$ Thus for
each fixed $y\in \Zp$ we obtain a formal power series $L(x,y)$ 
in $x^{-1}$ with
coefficients in a finite extension of $\KV$.
 
As mentioned above, it is now known that all characteristic $p$
$L$-series arising from
arithmetic have analytic continuations to ``essentially algebraic
entire functions'' on $S_\infty$ (see Subsection 8.5 of \cite{go1}).
That $L(s)$ is ``entire'' on $S_\infty$ means, in practice,
that for fixed $y\in \Zp$ every power series $L(x,y)$
is entire (i.e., converges for all values
of $x^{-1}$) {\it and} that the zeroes of this $1$-parameter family
of entire power series flow continuously.

The ``essential algebraicity'' rests on the following observation.
Let $\pi_\ast$ be a fixed $d_\infty$-th root of $\pi$ and let
$a\in \A$ be positive with $d(a)=d$. Let $j$ be a non-negative integer.
Notice that, by definition, 
$$(a)^{-(x\pi_\ast^j,-j)}=x^{-d}a^j\,;$$
that is, we have removed $\pi$ from the definition. Clearly the function
$x\mapsto (a)^{-(x\pi_\ast^j,-j)}$ is a polynomial in $x^{-1}$ with
algebraic coefficients. The essential algebraicity of $L(s)$ means that
the same thing happens with $L(s)$; i.e., the functions
$$z_L(x,-j):=L(x\pi^j_\ast,-j)$$
are polynomials (called the ``special polynomials'') 
with algebraic coefficients. 

Note that simple $p$-adic
continuity implies that $L(x,y)$ is the interpolation of $L(x,-j)=
z_L(x\pi^{-j}_\ast,-j)$ to arbitrary $y\in \Zp$. 

It is now a straightforward exercise to use the algebraic elements
$I^{-(\pi_\ast,-1)}$, for ideals $I$, to define $v$-adic analogs of the
local value field $\KV$ for finite primes $v$. We denote this $v$-adic
local value field by $\k_{v,\mathbf V}$.

Using the above notational conventions, we 
let $\mathfrak k$ be a finite field extension of
$k$ contained in our complete, algebraically closed extension
$C_\infty$ of $K=k_\infty$. Let $\phi$ be a Drinfeld module (or
an $\A$-module) that is defined over $\mathfrak{k}$.
Let $\mathfrak{p}$ be a finite prime of $k$ and let $\mathfrak{P}$ 
be a prime of $\mathfrak{k}$ lying over it with associated finite
field $\F_\mathfrak{P}$. Almost all such primes $\mathfrak{P}$
are ``good'' for $\phi$ in that reducing the coefficients of $\phi$ modulo
$\mathfrak{P}$ leads to a Drinfeld module $\phi^{(\mathfrak{P})}$
(or $\A$-module) over the finite
field $\F_\mathfrak{P}$
of the same ``rank'' as $\phi$ (i.e., the $\A$-ranks of the various
torsion modules prime to
$\mathfrak{P}$ are the same). In complete analogy with abelian varieties,
there is a ``Frobenius endomorphism,'' $Fr_\mathfrak{P}$,
of $\phi^{(\mathfrak{P})}$
and we set
$$f_\mathfrak{P}(t):=\det(1-tFr_\mathfrak{P} \mid T_v(\phi^{(\mathfrak{P})})),,$$
where $T_v(\phi^{(\mathfrak{P})})$ is the Tate module of
$\phi^{(\mathfrak{P})}$ at a prime $v\neq\mathfrak{p}$ of
$\A$. Again in analogy with abelian varieties, the polynomial
$f_\mathfrak{P}(t)$ has coefficients in $\A$ which are independent
of $v$ and has zeroes which satisfy the local Riemann hypothesis
(established by Drinfeld); see Subsection 4.12 of \cite{go1}. 
Let $n\mathfrak{P}\subseteq \A$ be the ideal norm of $\mathfrak{P}$; thus.
finally, one sets
$$L(\phi,s):=\prod_{\mathfrak{P}~\rm good} f_\mathfrak{P}(n\mathfrak{P}^{-s})^{-1}\,.$$
(See Part 3 of Remarks \ref{remarks1} for Euler factors at the finitely
many ``bad'' primes.)

As mentioned above, due to very recent work of G.\ Boeckle, R.\ Pink, 
and the author, it is known that
all $L(\phi,s)$ (and all associated partial $L$-series)
have an analytic continuation at
$\infty$ to an essentially algebraic
entire function. This analytic
continuation allows us to work in almost total generality in the function
field theory in obvious contrast to our need for abelian $L$-series
classically.

\begin{rem} \label{tagwan1}
In the basic case of $\A=\Fr[T]$, the analytic continuation of these
$L$-series was first obtained by Taguchi and Wan (\cite{tw1}, \cite{tw2}).
This was obtained by expressing the $L$-series as a Fredholm determinant
in the manner of classical Dwork theory. It is not known  at the present
time whether one can use these results in the manner of \cite{ks1}.
Wan's original paper on the characteristic $p$ Riemann hypothesis
\cite{w1}, which uses elementary methods, arose originally as a response
to a query from the present author as to whether the methods of
\cite{tw1}, \cite{tw2} could be used to obtain the exponential decay
of certain $L$-series coefficients which can be computed in closed
form (see Example 8.24.2 of \cite{go1}). In fact, as will be seen below,
this exponential decay occurs in complete generality.
\end{rem}

In classical theory, one constructs $p$-adic $L$-series out of the
special values of complex $L$-series. We will briefly describe here
how a very analogous construction works for the special polynomials
mentioned above; the reader may choose to ignore this construction 
during a first reading. In any case,
let $v$ be a finite prime of $\A$
and set $S_v:=\Z/(r^{\deg v}-1)\times \Zp$.
Let $(x,y)\in S_\infty$. Since the Euler factors in $L(\phi,s)$ have
coefficients in $\A$, it makes sense to also study them $v$-adically.
Thus, one defines the $v$-adic $L$-series
associated to $\phi$ etc., by using an Euler product over the finite primes of
$\mathfrak k$ {\bf not} lying over $v$ (in the obvious sense using the map
$\theta$). One obtains in this fashion an essentially algebraic
entire function $L_v(x_v,s_v)$ on the space $\bfC_v^\ast\times S_v$, where
$\bfC_v$ is the completion of a fixed algebraic closure of
$\k_v$ equipped with its canonical topology (see Subsection 8.3 of
\cite{go1}). However, it is also easy to see that this function is just
the $v$-adic interpolation of the special polynomials $z_L$ mentioned
above.  Thus we see how very close the theory is at all the places of
$\k$ and the central role played by the special polynomials.

\begin{example} \label{zetaA}
 Let $\A=\Fr[T]$, $k=\Fr(\theta)$. Note that ideals are in one to one
correspondence with monic polynomials. Thus we set
$$\zeta_A(s):=\sum_{n\in \A~\rm monic} n^{-s}=
\prod_{f~\rm{monic~prime}}(1-f^{-s})^{-1}\,,$$
If we expand out $n^{-s}$, we find
$$\zeta_A(x,y)=\sum_{n~\rm monic}x^{-\deg n}\langle n \rangle^{-y}=\sum_{j=0}^\infty x^{-j}\left(\sum_{\deg n=j}\langle n \rangle^{-y}\right)\,.$$
The function $\zeta_A(s)$ is obviously {\bf the} analog of the
Riemann zeta function for $\A$. It interpolates $v$-adically to
$$\zeta_{A,v}(x_v,s_v)=\prod_{f\neq v~\rm monic~prime}(1-x_v^{-\deg f}f^{-s_v})^{-1}=\sum_{j=0}^\infty
x_v^{-j}\left(\sum_{\substack{n{\rm ~prime~to~} v\\
\deg n=j}}n^{-s_v}\right)\,.$$
$\zeta_A(s)$ has special values and ``trivial zeroes'' in line with
what is known for the Riemann zeta function (see Subsections 8.12,
8.13 and 8.18
of \cite{go1}).
\end{example}

\begin{example} \label{Lcar}
Let $\A=\Fr[T]$ and let $C$ be the Carlitz module over
$k$ as in Example \ref{2carlitz}. The Tate modules of the Carlitz
module have rank $1$ and are very similar to the Tate modules of roots
of unity. Thus a simple calculation gives
$$L(C,s)=\zeta_A(s-1)\,.$$
There is an obvious $v$-adic version of this result.
\end{example}

Let $k^{\rm sep}\subset C_\infty$ be the separable closure of $k$
and let $G:={\rm Gal}(k^{\rm sep}/k)$. Let $\rho\colon G\to
GL_m(\bfC_\infty)$ be a representation of Galois type (i.e., factoring
through a finite extension of $k$). A completely similar theory holds for
the $L$-series $L(\rho,s)$ formed in the obvious fashion.

\begin{rems}\label{remarks1}
1. In practice, the essential algebraicity at $\infty$ 
of an $L$-series $L(s)$ 
given above actually arises as a consequence of its
being entire. Indeed, by construction $z_L(x-j)$ is certainly
an entire power series in $x^{-1}$. Now in most cases, such as the $L$-series
of a Drinfeld module, one concludes that $z_L(x,-j)$  also has
$\A$-coefficients. The only way this can happen is that almost all the
coefficients
are $0$; i.e., $z_L(x,-j)$ is a polynomial. The other cases usually
factor a function of this sort and so their essential algebraicity can
also be deduced.\\
2. The special polynomials $z_L(x,-j)$ lie at the heart of the analytic
continuation given by the author, G.\ Boeckle and R.\ Pink. Indeed, Boeckle
and Pink \cite{bp1}
give a cohomological expression for these polynomials which
leads to an estimate of the degree (in $x^{-1}$) of $z_L(x,-j)$. Moreover,
Boeckle has used this estimate to establish
that this degree grows {\it logarithmically} with $j$.
It is then relatively easy to translate this into a logarithmic
growth statement for the measures associated to $L(s)$ (at all places of
$\k$).
As explained in \cite{ya1}, this is enough to establish that the integrals
for the functions $L(s)$ (again, at all places of $\k$) converge everywhere.\\
3. It is obviously desirable to have Euler factors at all the finite
primes (as opposed to just the ``good'' primes). It seems likely that
the work of Boeckle and Pink will be able to provide these in line with
that is known for abelian varieties classically. In any case, it is 
easy to see that there are many examples of Drinfeld modules with
no bad primes, even when they are defined over $k$ (unlike the situation with abelian varieties).
\end{rems}

The proof of the
analytic continuation of the $L$-series $L(s)$ at the various places of
$\k$, mentioned in Part 2 of Remarks \ref{remarks1}, works quite similarly
whether the place is $\infty$ or a finite prime. This reinforces
previous experience that the theories at $\infty$ and at a finite prime $v$
are substantially the same. For instance, as mentioned
above they are all interpolations of the special polynomials $z_L(x,-j)$.

Therefore it is somewhat reasonable to expect that an analog of
Conjecture \ref{classGRH1} should work for {\it all} the
interpolations of $L(s)$ at all the places of $\k$.
As we have mentioned when $r=p$, Daqing Wan calculated 
in \cite{w1} the Newton polygons
at $\infty$ of the power-series $\zeta_A(x,y)$ (as in
Example \ref{zetaA}) for all $y\in \Zp$. Wan
found that these polygons were always  simple implying that
the zeroes (in $x^{-1}$) of $\zeta_A(x,y)$ are always
in $\K=\k_\infty$ and are themselves  simple.  As $\R=\Q_\infty$, 
this is obviously in keeping with the
classical Conjectures \ref{classGRH2} and \ref{clsimpl}.

Wan's proof was simplified by D.\ Thakur and J.\ Diaz-Vargas in \cite{dv1}.
Finally, based on some work of B.\ Poonen for $r=4$, the general
case (all $r$) was established by J.\ Sheats in \cite{sh1}. Still, it was
not clear how to proceed to a ``good'' version of the
Generalized Riemann Hypothesis, etc., in the function field case. For instance,
there are examples of zeroes of
$\zeta_A(s)$ (in the obvious definition) which
do not belong to $\KV$ for some non-polynomial
$\A$.
Moreover, Daqing Wan had mentioned to the author that the calculation
at $\infty$ for $\zeta_{\Fr[\theta]}(s)$, could 
easily be modified to establish the {\it same}
result $v$-adically when $\deg v=1$ (or at least for those $s_v$ in
an open subset of $S_v$ for such $v$; see
Proposition \ref{wans}).

It is well-known (and important for what follows in the next
section) that
non-Archimedean analysis is quite algebraic when compared to
complex analysis. This is
brought out by the result that {\it all} entire functions are
determined up to a constant by their zeroes (with multiplicity of course)
and that these zeroes are algebraic over any complete field containing
the coefficients. In particular it makes sense to analyze such functions
by studying the extension fields obtained by adjoining the zeroes.
Let $L(s)$,
$s=(x,y)\in S_\infty$, be one of the characteristic $p$ $L$-series 
arising from arithmetic (via a Drinfeld module or $\A$-module)
and let $\K_L$ be the finite
extension of $\KV$ obtained by adjoining the coefficients of $L(s)$. (In 
practice, $\K_L$ will usually be $\KV$ itself, some finite constant
field extension of $\KV$ obtained by adjoining the values of 
certain characters, or some finite extension obtained by adjoining
``complex multiplications'' etc.) 
Let $\K_L(y)$ be the extension of $\K_L$ obtained
by adjoining the roots of $L(x,y)$ for each $y$. Thus, a-priori,
$\K_L(y)$ is merely some algebraic extension of $\K_L$.

Now in any extension of local, or global, function fields, the most
important part is the maximal separable subfield. Indeed it is
well known that the total extension is uniquely determined by
its degree over the maximal separable subfield 
(see, eg., 8.2.12 of \cite{go1}). Thus a first function field version
of Conjecture \ref{classGRH1} (or Conjecture \ref{classGRH2}) is the
following.

\begin{conj} \label{ffieldGRH1}
The maximal separable (over $\K_L$) sub-extension of $\K_L(y)$ is
{\it finite} over $\K$ for all $y\in \Zp$. Similarly, the maximal
separable subfields of the extensions obtained by adjoining the
$v$-adic zeroes should also be finite for each $s_v\in S_v$.
\end{conj}

\begin{rems}\label{notapply}
1. Notice that the conjecture is vacuously true for
$y=-j$ since, obviously, the special polynomials have only finitely
many zeroes to begin with.\\
2. At $\infty$ one can show that the maximal separable (over
$\K$) subfield of $\K_L(y)$ is independent of the choice of sign function
or uniformizer. Thus it depends  only on $y$ and the underlying ``motive''
used to construct the $L$-series. (The formal module of such a motive may
always be extended to separable elements. We believe that these extended
modules may be of importance in studying the zeroes of characteristic $p$
$L$-series.)
\end{rems}

Conjecture \ref{ffieldGRH1} is clearly in line with the few known results
about such extensions (e.g., the results of Wan, Diaz-Vargas etc.). However
it suffers from the drawback of not being a statement directly
about the zeroes themselves. On the other hand, the analytic
continuation of these functions using integral calculus (mentioned above)
is based on the a-priori estimates of Y.\ Amice in
non-Archimedean functional analysis; these estimates are nicely reviewed in
\cite{ya1}. In particular, the logarithmic
growth of the degrees of the special polynomials translates into the
exponential decay of the coefficients of the $L$-series when expressed
as power-series as above. (In the particular case 
where $\A=\Fr[T]$ such exponential
decay for $L$-series of Drinfeld modules, etc., was originally
mentioned to the author by D.\ Wan.)

\begin{rem} \label{consexpdecay}
The exponential decay of the coefficients
has a number of consequences. For instance, let
$a\in \A$ be non-constant. One can  construct the ``Carlitz polynomials''
(as in \cite{ca1} or Subsection 8.22 of \cite{go1})
for the ring $\Fr[a]\simeq \Fr[T]\subseteq \A$ by simple substitution of
$a$ for $T$. Then the above mentioned
exponential decay and the main result of \cite{ca1} imply that all
$L$-series arising from arithmetic have convergent expansions in these
Carlitz polynomials at $\infty$. Such expansions are somewhat similar
to Fourier expansions classically.
\end{rem}

The exponential decay also 
suggests refining Conjecture \ref{ffieldGRH1} in the following fashion.
Let us write
$$L(x,y)=\prod_i\left(1-{\beta_i^{(y)}}/{x}\right)\,,$$
where the elements $\{\beta_i^{(y)}\}$ are algebraic over
$\K_L$ and are the zeroes of $L(x,y)$. Clearly we are only interested
in non-zero $\beta_i^{(y)}$ and we call
$\{\lambda_i^{(y)}:=1/\beta_i^{(y)}\}$ the ``reciprocal zeroes''
of $L(s)$; standard non-Archimedean function
theory tells us that the reciprocal zeroes
are discrete in the sense that their absolute values
tend to infinity. There is an obvious
$v$-adic version of this product. We are thus lead to the following
refinement of Conjecture \ref{ffieldGRH1} which focuses directly
on the zeroes.

\begin{conj}\label{ffieldGRH2}
There exists a positive real number $b=b(y)$ such that
if $\delta\geq b$, then there is at most {\bf one} reciprocal
zero (taken without multiplicities) 
of $L(x,y)$ of absolute value $\delta$. (This conjecture,
again, is meant to apply to all the interpolations of $L(s)$
at all the places of $\k$.)
\end{conj}
 
\noindent
In other words, for fixed $y\in \Zp$ almost all zeroes are uniquely
determined by their absolute values etc.

\begin{rems}\label{caveats1} 
1. Conjecture \ref{ffieldGRH2}
is also vacuously true for $y=-j$.
Perhaps a more refined version might be non-trivial in this case also.\\
2. Conjecture \ref{ffieldGRH2} appears to us to be the correct version
of the Generalized Riemann Hypothesis in the characteristic $p$ setting.
However, it should probably still be viewed as a working version since we
do not yet even know what the implications of this conjecture are
(unlike, obviously, classical theory).\\
3. Let $f(t)$ be an entire power series with coefficients in a finite
extension of $\Qp$. As there are only finitely many extensions of $\Qp$ of
bounded degree, {\it any} finite bound on the number of zeroes of $f(t)$ of
fixed degree will suffice to establish that the zeroes of $f(t)$ generate
a finite extension of $\Qp$. This argument does not work in finite
characteristic however. As the reader will see, Conjecture \ref{ffieldGRH2}
is the best one can hope for in finite characteristic;
see Example \ref{sqrtcar} of our next section.
\end{rems}

Our next result establishes that Conjecture \ref{ffieldGRH2} implies
Conjecture \ref{ffieldGRH1}. We show it at $\infty$ with the $v$-adic
result being completely analogous.

\begin{prop} \label{ffGRH2toffGRH1}
Assume {\rm Conjecture \ref{ffieldGRH2}} and let $b(y)$ be as in its statement.
Let $\beta_i^{(y)}$ be a root such that $|\lambda_i^{(y)}|=
1/|\beta_i^{(y)}|>b(y)$. Then $\beta_i^{(y)}$ is totally inseparable
over $\K_L$.
\end{prop}

\begin{proof} Let $\overline{\K}\subset \bfC_\infty$ be the algebraic
closure. Let $\sigma$ be an automorphism of $\overline{\K}$ which
fixes $\K_L$. Then $\sigma(\beta_i^{(y)})$ is also a zero
of $L(x,y)$ of the {\it same} absolute value. Thus Conjecture
\ref{ffieldGRH2} implies that $\sigma(\beta_i^{(y)})$ actually equals
$\beta_i^{(y)}$. This is equivalent to the total inseparability of
$\beta_i^{(y)}$.
\end{proof}

\noindent
Thus we only have finitely many zeroes that might contribute separable
elements and hence
the maximal separable subfield is obviously finite. This is 
Conjecture \ref{ffieldGRH1}. Ultimately, one would like an arithmetic
characterization of this maximal separable subfield.

It is explicitly allowed in Conjecture \ref{ffieldGRH2} that one can
throw out finitely many zeroes. This is necessary for the Conjecture
to apply at all the places of $\k$. Indeed, when one interpolates
$v$-adically, one removes the Euler-factors lying over $v$. This process
in fact adds finitely many zeroes $v$-adically which may be arbitrary
in behavior, as our next example attests.

\begin{example} \label{vadiczeroes}
Let $\zeta_A(s)$ be as in Example \ref{zetaA} which we write
as
$$\zeta_A(s)=\sum_{n~\rm monic}n^{-s}\,.$$
So, for $j\geq 0$ we have
$$z_\zeta(x,-j)=\sum_{j=0}^\infty x^{-j}(\sum_{\deg n=j}n^j)\,.$$
As mentioned, these power-series are actually polynomials in
$x^{-1}$. Let $v$ be a finite prime associated to an irreducible monic
$f$ of degree $d$. Let $s_v\in S_v$. To interpolate $v$-adically, one takes a sequence $e_t$
of natural numbers with the property that $e_t$ goes to $\infty$ in
the Archimedean absolute value but to $s_v$ in the $p$-adic absolute value.
In particular, this process will eliminate all monic $n$ which
are divisible by $f$. Thus $z_\zeta(x,-j)$ is transformed $v$-adically into
$$(1-x_v^{-d}f^j)z_\zeta(x_v,-j)\,.$$
Therefore we have added the zeroes of $1-x_v^{-d}f^j$ to the
zeroes of $z_\zeta(x_v,-j)$. Note that there will be many zeroes of
$1-x_v^{-d}f^j$ of the same size if $d>1$ is prime to $p$.
\end{example}

Example \ref{vadiczeroes} also makes clear the importance of having Euler
factors at {\it all} finite primes in the definition of $L$-series
and not just the good primes; indeed,
we add (finitely many) zeroes for each Euler factor left out. A refined
version of these conjectures should ultimately take this into account.
In any case,
throwing out finitely many zeroes also allows
Conjecture \ref{ffieldGRH2} to have some surprising consequences that make
the characteristic $p$ and classical cases seem even closer; see Section
\ref{miss}.

Conjecture \ref{ffieldGRH2} has a direct analog
in classical theory.

\begin{conj} \label{classGRH3}
Let $e\geq 0$. Then there are at most {\bf two} zeroes (taken without
multiplicity) of $\Xi(\chi,t)$ of absolute value $e$.
\end{conj}

\begin{prop} \label{uptofinitelymany}
{\rm Conjecture} \ref{classGRH3} implies the Generalized Riemann
Hypothesis up to possibly finitely many exceptional
zeroes (which are then the non-critical real
zeroes of $\Lambda(\chi,s)$).
\end{prop}

\begin{proof}
We begin by assuming that $\chi$ is a real valued character.
Let $z$ be a zero of $\Xi(\chi,t)$. Then
complex conjugation, and the functional equation, imply
that $z$, $-z$, $\bar{z}$ and $-\bar{z}$ are also zeroes of $\Xi
(\chi,t)$ with of course the same
absolute value. If $z$ is non-real and non-purely imaginary, 
then these elements are distinct
which is a violation of the conjecture. The purely imaginary
zeros of $\Xi(\chi,t)$ are the real zeroes of $\Lambda(\chi,s)$ of
which there are only finitely many by the functional
equation. 

We postpone the proof when $\chi$ is non-real to the next section.
\end{proof}

\begin{rems} \label{feuse}1.
The reason that Conjecture \ref{classGRH3} allows two zeroes where
Conjecture \ref{ffieldGRH2} allows only one lies with 
functional equations. The functional equation
of $L(\chi,s)$ precisely allows us to
move the critical line over to the real axis, and so make the connection to
the characteristic $p$ theory. On the other hand, it also introduces
a new symmetry for $L(\chi,s)$ and so adds to the count of zeroes.\\
2. The use of absolute values in these conjectures offers a unification
between the global and local Riemann Hypotheses (by ``local'' we mean a
Riemann Hypothesis for a variety or Drinfeld module etc., over a finite field).
Indeed,  all local Riemann Hypotheses are stated in terms of absolute
values of the zeroes.\\
3. Let $\chi$ be a real-valued character associated to a number field. 
It is easy to see that, modulo
possible real non-critical zeroes of $L(\chi,s)$,
Conjecture \ref{classGRH3} is equivalent to Conjecture \ref{classGRH2}.
However, in the function field case, it is also easy to see that
Conjecture \ref{ffieldGRH2} is much stronger than merely assuming that almost
all zeroes are totally inseparable over $\K_L$ etc. \\
4. It is known that there are no real zeroes for the Riemann zeta function.
In any case, one can certainly use the sign of the real part of a zero
to pose a refinement of Conjecture \ref{classGRH3} which will work for
all zeroes. Only time will tell if there is any utility in such
statements.\\
5. As in the proof of Proposition \ref{uptofinitelymany}, Conjecture \ref{classGRH3}
may be restated as saying that a zero of $\Xi(\chi,t)$ is uniquely
determined up to multiplication by $\pm 1$ by its absolute value.
\end{rems}

Crucial to the proof of Proposition \ref{uptofinitelymany} is the fact
that functional equation of $\Xi(\chi,t)$ takes $t$ to $-t$ and so
is {\it absolute value preserving}. This suggests very strongly that
any sort of functional equation for our characteristic $p$ $L$-series should
also be given by invariance under 
mappings of this sort. However, any non-trivial such automorphism would
automatically provide more than $1$ zero of a given absolute value. 
Perversely, therefore, by Part 3 of Remarks \ref{caveats1} we see that
a classical style functional equation in characteristic $p$ would lead
to a Riemann hypothesis with {\it no} punch!

It remains to present the analog of the Generalized Simplicity Conjecture in
the characteristic $p$ case. We first make some definitions.
Let $\overline{k}$ be a fixed algebraic closure of
$k$ and let $k^{\rm sep}$ be the separable closure. Let
$G:={\rm Gal}(k^{\rm sep}/k)$ and let $\rho\colon G \to \bfC_\infty^\ast$
be a character which factors through a finite abelian extension. (We
restrict ourselves to abelian characters for comparison with the classical
case where we are only using abelian characters; however see below). 
Let $G_1$ be the Galois group of the maximal constant field extension
of $\K$ over $\K$ and let $\sigma\in G_1$. Obviously
$\rho^\sigma$ is also a character of $G$ of the same type.
We call the orbit of $\rho$ under these automorphisms the
{\it Galois packet} generated by $\rho$.

The Galois packet of a classical character $\chi$ is just $\{\chi,
\bar{\chi}\}$.

The characteristic $p$ version of the Simplicity Conjecture can now
be given. We state it at $\infty$ with the obvious
$v$-adic version being left to the reader.

\begin{conj} \label{charpsim}
1. Let $L(s)=L(\rho,s)$, $s=(x,y)$, be the characteristic $p$ $L$-series
of $\rho$. Then for fixed $y$, almost all zeroes of $L(x,y)$ are simple.\\
2. Let $\rho$ and $\rho^\prime$ be two characters as above but from
distinct Galois packets. Let $y\in \Zp$. Then there is a
number $c=c(y)$ such that the absolute values of the zeroes of
$L(\rho,x,y)$ and $L(\rho^\prime,x,y)$ which are $>c$ are distinct.
\end{conj}

\begin{rem} \label{simplemotive}
The abelian representation $\rho$ may be viewed as a ``simple motive over
$k$;'' that is, $L(\rho,s)$ has no ``obvious'' factors. (Indeed,
as the degree of $\rho$ is $1$, every Euler factor is linear and so
obviously does not factor.) Similarly, due to the work of Richard
Pink \cite{p1} on the Serre Conjecture for Drinfeld modules, a Drinfeld 
$\Fr[T]$-module over $\Fr(\theta)$ is also a simple motive.

Now let $\overline{\Qp}$ be an algebraic closure of $\Qp$ and let
$V$ be a finite dimensional $\overline{\Qp}$-vector space. Let $G=
{\rm Gal}(k^{\rm sep}/k)$ and let $\rho$ now be an irreducible representation
of $G$ into ${\rm Aut}_{\overline{\Qp}}(V)$ of Galois type (i.e., which
factors through a finite Galois extension of $k$). As explained
in Subsection 8.10 of \cite{go1}, one works with such characteristic
$0$ representations in order to get the correct 
(=classical) local factors at the
ramified primes. Let $R_p$ be the integers of $\overline{\Qp}$ and
let $M_p$ be the maximal ideal. Since $\rho$ factors through a Galois
representation, the coefficients of the classical Euler factors at
any prime obviously belong to $R_p$. The corresponding {\it characteristic
$p$} $L$-series is then defined by reducing the classical
Euler factors modulo $M_p$ and viewing
the reduced polynomial as having $\mathbf{C}_\infty$-coefficients; one then
forms the obvious Euler product. In
a similar way we can reduce $\rho$ itself modulo $M_p$. If this {\it reduced}
representation is still irreducible, then $\rho$ gives rise to
a simple (characteristic $p$!) motive and the Simplicity Conjecture
should also be true for such $\rho$.
\end{rem}

We shall discuss Conjecture \ref{charpsim} in more detail at the end of
our next section (after Corollary \ref{c4}).
However, note that even classically one throws out finitely many
zeroes in the Generalized Simplicity Conjecture (e.g., $s=1/2$)!

\section{Missing Zeroes and the Generalized Simplicity Conjecture} \label{miss}

Let $p(z)$ be a polynomial with coefficients in
some field and with $p(0)=1$. Basic high school algebra tells us
that every field containing the roots of $p(z)$ must also contain
its coefficients. Now let $p(z)$ be an entire non-Archimedean function
with $p(0)=1$ (all $L$-series in the characteristic $p$ theory are of this
form). As $p(z)$ factors over its zeroes, we see analogously that
every {\it complete} field containing the zeroes also contains the
coefficients.

As is universally known the above discussion
is totally false for complex entire functions.
Indeed, complex
entire functions need not have zeroes at all (e.g., $e^z$) or the zeroes
can lie in strictly smaller complete fields than the coefficients
(e.g. $1-e^{2\pi i z}$). From the point of view of non-Archimedean
analysis such entire functions have ``missing zeroes.''

In this section we will show how Conjecture \ref{ffieldGRH2} implies
an analog of ``missing zeroes'' for the characteristic $p$ $L$-functions.
That is, we will show how it implies that almost all zeroes lie in
strictly smaller subfields than one would a-priori believe they
should -- this makes the function field theory and the classical
theory look even closer than one would at first believe. 
In the process we will re-examine the Simplicity Conjectures for
both function fields and number fields. We will do this by examining similar
abelian
$L$-series in tandem for both cases. The reader will then see how to
handle other $L$-functions in both theories.

As above, let $G={\rm Gal}(k^{\rm sep}/k)$ and 
Let $\rho\colon G\to \bfC_\infty^\ast$ be
a character arising from a finite abelian extension
$\mathfrak{ k}=\mathfrak{k}_\rho$ of 
$k$. We explicitly assume that the values of
$\rho$ are {\bf not} contained in $\K$. Thus $\rho$ is analogous
to a non-real classical abelian character $\chi$.
We fix one such complex character $\chi$ and let $\mathfrak {L}
/\mathfrak{L}_1$ be a
finite abelian extension of number fields
such that $\chi$ is defined on ${\rm Gal}(\mathfrak{L}/
\mathfrak{L}_1)$. We do not
assume
that the base $\mathfrak{L}_1$ of this extension is $\Q$ (while, for reasons which will
be clear later, we {\bf do} make the analogous assumption in the function
field case). We let $L(\rho,s)$, 
$s\in S_\infty$, be the $L$-series of $\rho$ at $\infty$ (which will be
enough for our purposes) and let $L(\chi,s)$, $s$ a complex number, be the
classical $L$-series of $\chi$. 
 
The symbol ``$\chi$'' will always be reserved for complex abelian
characters of number fields and ``$\rho$'' will be reserved for
characteristic $p$ valued characters of function fields.

We begin by applying Conjecture \ref{ffieldGRH2} directly to
$L(\rho,s)$ where $s=(x,y)\in S_\infty$. Notice first that
$\K_L$ clearly equals $\KV(\rho):=$ the constant field extension of $\KV$ obtained
by adjoining the values of $\rho$.
As in Proposition \ref{ffGRH2toffGRH1}, we deduce immediately 
the following result.

\begin{prop} \label{c1}
Let $y\in \Zp$ be fixed. Then 
{\rm Conjecture \ref{ffieldGRH2}} implies  that
almost all zeroes of $L(\rho,x,y)$ are totally inseparable over 
$\K_L$.
\end{prop}

Now we examine $L(\chi,s)$ or rather
$\Xi(\chi,t)$. Note firstly that the analog of $\KV$ in this
case is obviously
$\R$ and the analog of $\K_L$  is obviously $\C$ (which equals 
$\R$ adjoined with the values
of $\chi$). The functional equation relates $\Xi(\chi,t)$
and $\Xi(\bar{\chi},-t)$. Suppose that $z$ is a zero of $\Xi(\chi,t)$.
Playing off the functional equation for $\Xi(\chi,t)$ and complex conjugation,
we deduce that $\bar{z}$ is also a zero for $\Xi(\chi,t)$. This is all
the information that one may deduce directly. Thus applying
Conjecture \ref{classGRH3}  directly
to $\Xi(\chi,t)$ implies that the zeroes
of $\Xi(\chi,t)$ belong to $\C$; i.e., nothing interesting!

Now let us return to the characteristic $p$ theory. We introduce
another ingredient into the mix by simply noting that $L(\rho,s)$
 divides the zeta function of the ring of integers
$\mathcal{O}_\mathfrak k$ of $\mathfrak k$. This zeta function is 
also an essentially algebraic entire function. Note that $\K_\zeta=\KV$
which is strictly smaller than $\K_L$.
 By applying  Conjecture
\ref{ffieldGRH2} to $\zeta_{\mathcal{O}_\mathfrak k}(s)$
we deduce immediately the next result which improves dramatically Proposition
\ref{c1}.
\begin{prop} \label{c2}
Let $y\in \Zp$ be fixed. Then 
{\rm Conjecture \ref{ffieldGRH2}} implies  that
almost all zeroes of $L(\rho,x,y)$ are totally inseparable over 
$\KV$.
\end{prop}
\begin{proof}
By Proposition \ref{c1},
Conjecture \ref{ffieldGRH2} implies that almost all zeroes of
 $\zeta_{\mathcal{O}_\mathfrak k}(x,y)$ are totally
inseparable over $\KV=\K_\zeta$. This is then obviously true for
$L(\rho,x,y)$.
\end{proof}

\begin{rems}\label{evenbetter}
1. In fact, we shall do even better once we bring the first part of the 
Generalized Simplicity
Conjecture into the mix; see Proposition \ref{c6}.\\
2. Notice that, as $\KV$ is totally inseparable over $\K$, Proposition
\ref{c2} immediately implies the total inseparability of almost
every zero all the way down to
$\K$ itself.\\
3. Again we remind the reader that there are $v$-adic versions to
the above results.
\end{rems}

Next, let us return to the classical case. We have seen that
for $\Xi(\chi,t)$, Conjecture \ref{classGRH3} says nothing interesting
directly. However, we shall now play the same game as we did in
characteristic $p$ and we shall obtain an analogous improvement. 
We know that our $L$-series divides the
zeta function $\zeta_\mathfrak{L}(s)$ of $\mathfrak{L}$. We shall denote
the ``$\Xi$-version'' of this zeta function by $\Xi(\chi_0,t)$ where
$\chi_0$ is the appropriate trivial character.

\begin{prop}\label{c3}
{\rm Conjecture \ref{classGRH3}} implies that
almost all zeroes of $\Xi(\chi,t)$ are real.
\end{prop}

\begin{proof}
As above we deduce that 
Conjecture \ref{classGRH3} implies that almost all zeroes
of $\Xi(\chi_0,t)$ are real. We therefore immediately deduce
the same for $\Xi(\chi,t)$.
\end{proof}

\noindent
From Part 3 of Remarks \ref{feuse} we obtain the following corollary.

\begin{cor} \label{c4} Let $\chi$ be an abelian character of
a number field. Then, modulo possible real non-critical zeroes of $L(\chi,s)$
(of which there can be only finitely many),
{\rm Conjecture \ref{classGRH3}} is equivalent to the GRH for
$L(\chi,s)$.
\end{cor}

There is even more that can be gleaned from this line of thought.
Let $\rho$, $\rho^\prime$ be two characteristic $p$ abelian characters as
above.

\begin{prop} \label{c5}
Assume that $\rho$ and $\rho^\prime$ are in the same Galois packet
and fix $y\in \Zp$. Then {\rm Conjecture \ref{ffieldGRH2}} implies that
almost all zeroes of $L(\rho,x,y)$ and $L(\rho^\prime, x,y)$
are equal.
\end{prop}

\begin{proof}
Let $\sigma$ be the automorphism taking $\rho$ to $\rho^\prime$.
Extend $\sigma$ to the full algebraic closure. Let $\beta$ be a zero
of $L(\rho,x,y)$ where $1/|\beta|$ is sufficiently large so
that Conjecture \ref{ffieldGRH2} applies for both
$L(\rho,x,y)$ and $L(\rho^\prime,x,y)$. Then $\sigma (\beta)=\beta$ by
Proposition \ref{c2}. On the other hand, $\sigma (\beta)$ is a zero
of $L(\rho^\prime,x,y)$. Thus by Conjecture \ref{ffieldGRH2} it
must be the unique zero of its absolute value.
\end{proof}

\noindent
The obvious $v$-adic version of Proposition \ref{c5} also holds.

Proposition \ref{c5} explains why Part 2 of the function field 
Generalized Simplicity
Conjecture (Conjecture \ref{charpsim}) is formulated as it is. 

\begin{rems} \label{ruleout}
1. Proposition \ref{c5} rules out a general elementary proof by Newton
Polygons of Conjecture \ref{ffieldGRH2}.\\
2. Proposition \ref{c5} therefore gives some credence to viewing the
Galois action as being the characteristic $p$ analog of {\it both}
the classical action of complex conjugation and the functional equation.
\end{rems}

Our next result shows how close the classical and function field
Generalized Simplicity Conjectures are.

\begin{prop} \label{c5.1}
Let $\chi$ and $\chi^\prime$ be two complex characters from different
Galois packets (i.e., $\chi$ and $\chi^\prime$ are not complex conjugates).
Then the Generalized Simplicity Conjecture ({\rm Conjecture
\ref{clsimpl}}) and the
Generalized Riemann Hypothesis imply that all zeroes, except possibly $t=0$,
of $\Xi(\chi,t)$ and $\Xi(\chi^\prime,t)$ have distinct absolute values.
\end{prop}

\begin{proof}
Assume GRH. Then, using complex conjugation and
the functional equation as before, we see that the $\Xi(\chi,t)$ and
$\Xi(\bar{\chi},t)$ fill out all zeroes of a fixed absolute value.
Thus the result follows immediately from GSC.
\end{proof}

\noindent
Thus the difference between the classical and function field Generalized
Simplicity Conjectures lies in comparing $L$-series associated to
characters in the same Galois packet. In characteristic zero we expect
infinitely many
different zeroes, whereas by Proposition \ref{c5}, in characteristic $p$
we expect only finitely many different zeroes.

The function field theory appears to offer insight
into the arithmetic meaning of simplicity of zeroes. Indeed, let us return
to the case of $L(\rho,s)$ as above. Recall that we have supposed
that $\rho \colon {\rm Gal}(k^{\rm sep}/k)\to 
\bfC_\infty^\ast$. Applying the function field Generalized 
Simplicity Conjecture now gives the next result.

\begin{prop} \label{c6}
The first part of {\rm Conjecture
\ref{charpsim}} (assumed along with
{\rm Conjecture \ref{ffieldGRH2}})  implies that, for fixed $y\in \Zp$ almost
all the zeroes of $L(\rho,x,y)$  belong to $\KV$.
\end{prop}

\begin{proof}
The first part of Conjecture \ref{charpsim}
(applied to $L(\rho,s)$) implies that almost all the
zeroes are simple. Combining this with
Conjecture \ref{ffieldGRH2}, and a simple argument using Newton Polygons,
we find that
 almost all zeroes of $L(\rho,x,y)$ are actually in $\K_L$.
On the other hand these zeroes are also totally inseparable
over the subfield $\KV$ of $\K_L$ by Proposition \ref{c2}. Note that
$\K_L=\KV(\rho)$ is separable over $\KV$. Thus, the only
way that this can happen is that the zeroes belong to $\KV$.
\end{proof}

\noindent
The obvious $v$-adic version of the proposition is also true.
As $\KV$ is strictly
smaller than $\K_L$ we have deduced an analog of the missing zeroes
phenomenon for function fields! Other examples can be easily worked out.
While it is 
true that in non-Archimedean analysis
we cannot have {\it all} zeroes in too small a field, there is
in fact no contradiction. Indeed, we 
are allowed in our conjectures to throw  out finitely many zeroes which avoids any difficulties
(as simple examples attest).

To finish, we will present an example which involves
complex multiplication by totally inseparable elements. This
leads to certain $L$-series where the $v$-adic versions of
Conjectures \ref{ffieldGRH2} and \ref{charpsim} are true,
but where the obvious analog of the $v$-adic version of Proposition \ref{c6} is not valid. 
That is, where one can find infinitely many non-trivial (over the 
$v$-adic analog of $\KV$)
totally inseparable roots for a fixed $s_v\in S_v$.
This highlights the crucial
role that separability plays in Proposition \ref{c6}.
In order to do so,
we first present Wan's cogent observation that the $\infty$-adic techniques 
presented in \cite{w1}, \cite{dv1}, and \cite{sh1} for $\Fr[T]$ also work
$v$-adically when $\deg v=1$.

\begin{prop} \label{wans}
Let $v$ be prime of degree $1$ in $\A=\Fr[T]$. Let $j$ be an element of the 
ideal $(r-1)S_v$. Then the zeroes (in $x_v$) of $\zeta_{A,v}(x_v, j)$ are
simple, in $\k_v$ and uniquely determined by their absolute value.
\end{prop}

\begin{proof} As $\deg v=1$, it is clear that $\K$ and $\k_v$ are isomorphic.
Moreover, without loss of generality,
we can set $v=(T)$. As the Newton polygon does not
depend on the choice of uniformizer, we choose our positive uniformizer
to be $\pi=1/T$ and we begin by letting $j$ be a positive integer divisible
by $r-1$.
Now the coefficient of $x^{-d}$ in $\zeta_A(x,-j)$ is precisely the
sum of $\langle n\rangle^j$ where $\deg n=j$ and $n$ is monic. On the other
hand, the
coefficient of $x_v^{-d}$ in $\zeta_{A,v}(x_v,-j)$ is the sum of
$n^j$ such that $n$ is monic of degree $d$ and $n\not \equiv 0~{\rm mod}\,v$.
This last condition is the same as saying that $n$ has non-vanishing
constant term. 

The set $\{\langle n \rangle\}$, where $n$ is monic, ranges over
all polynomials $f(1/T)$ in $1/T$ with constant term $1$ and degree (in $1/T$) 
$<d$. Moreover, as $j$ is divisible by $r-1$, the set
$\{\langle n \rangle^j\}$ is the same as the set
$f(1/T)^j$ where $f(u)$ is a monic 
polynomial of degree $<d$ and has non-vanishing
constant term. 

Let us denote by $\zeta_{A,v}(x,-j)$ the function obtained by replacing
$x_v$ by $x$ in $\zeta_{A,v}(x_v,-j)$ and applying the isomorphism
$\k_v \to \K$ given by $T\mapsto 1/T$. The above now implies
that
$$(1-x^{-1})^{-1}\zeta_{A,v}(x,-j)=\zeta_A(x,-j)\,.$$
The result for positive $j$ divisible by $r-1$ follows immediately.
The general result then follows by passing to the limit.\end{proof}

For the example, we let $\A=\Fr[T]$ with $r=2$. Let
$\A^\prime: =\Fr[\sqrt{T}]$, $\k^\prime:=\Fr(\sqrt{T})$ and 
$\K^\prime:=\Fr((1/\sqrt{T}))$. Form $A^\prime$, $k^\prime$ etc., in the
obvious fashion. Note that in this case $\KV=\K$. 
For each $a\in \A$ let $a^\prime\in \A^\prime$ be its unique
square-root. Let $\pi$ be a positive
uniformizer for $\K$ and let $\pi^\prime=\sqrt{\pi}$ be the uniformizer
in $\K^\prime$. Finally let $v=(g)$ be a prime of $\A$ of degree $1$ with
$v^\prime=(g^\prime)$ the unique prime of $\A^\prime$ above it.

\begin{example}\label{sqrtcar}
Let $\psi$ be the Drinfeld $\A$-module defined over $k^\prime$
given by
$$\psi_T(\tau):=\theta\tau^0+(\theta+\sqrt{\theta})\tau+
\tau^2\,.$$
It is simple to check that $\A^\prime$ acts as complex multiplications
of $\psi$; indeed, $\psi$ is just the Carlitz module $C^\prime$ for
$\A^\prime$ ($C^\prime_{\sqrt{T}}(\tau)=\sqrt{\theta}\tau^0+\tau$) as one readily checks.  Let $L(\psi,s)$, $s\in S_\infty$ be the
$L$-series of $\psi$ over $k^\prime$. As $\psi$ has complex multiplication,
$L(\psi,s)$ factors into the product of $L$-series associated to Hecke
characters. In this case, it is simple to work out what happens directly.
For each monic prime $f\in \A$, let $f^\prime$ be its unique
square-root in $\A^\prime$. Let 
$$L(s):=\prod_{f~\rm monic ~prime}(1-f^\prime f^{-s})^{-1}\,.$$
It is easy to see that $L(s)$ is the $L$-series of the Hecke character
$\Theta$ with $\Theta(f)=f^\prime$.
One checks readily that
$$L(\psi,s)=L(s)^2\,.$$
Thus we need only focus on $L(s)$. 

It clear that $\K_L$, for our $L$-series $L(s)$, equals $\K^\prime$
(=$\K(\Theta)$ defined in the obvious fashion) and so
is totally inseparable over $\K=\KV$.  Upon expanding 
$L(s)$ we find
$$L(s)=\sum_{n~\rm monic} n^\prime n^{-s}\,.$$
 Thus,
$$L(s)=\sum_n x^{-\deg n}n^\prime \langle n\rangle^{-y}\,.$$
Let $\langle z\rangle_{\pi^\prime}$ be the $1$-unit part of an element
of $z\in \K^\prime$ defined with respect to $\pi^\prime$.
Thus 
$$L(s)=\sum_{n^\prime ~\rm monic} x^{-\deg n} n^\prime \langle n^\prime \rangle_{\pi^\prime} ^{-2y}\,.$$
But notice that the degree in $T$ of $n$ is clearly the degree
in $\sqrt{T}$ of $n^\prime$. Thus, finally,
$$L(s)=\sum_{n^\prime ~\rm monic} x^{-\deg n^\prime} n^\prime 
\langle n^\prime \rangle_{\pi^\prime}^{-2y}\,.$$
If we form $\zeta_{A^\prime}(s)$ in the obvious fashion, then we have
shown that
$$L(s)=\zeta_{A^\prime}(x\pi^\prime, 2y-1)\,.$$
Thus the results of Wan, Thakur, and Diaz-Vargas tell us that
the zeroes of $L(s)$ are simple, in $\K^\prime=\K(\Theta)$ and
uniquely determined by their absolute value. As such they
are indeed totally inseparable over $\K$. Thus both Conjecture
\ref{ffieldGRH2} and Conjecture \ref{charpsim} are true for $L(s)$.

We now form the $v$-adic functions $L_v(\psi,x_v,s_v)$, $L_v(x_v,s_v)$,
etc. Note that, by definition,
$$L_v(x_v,s_v)=\prod_{f\neq g}(1-f^\prime x_v^{-\deg f}f^{-s_v})^{-1}\,.$$
Note further that the $v$-adic version, $\k_{v,\mathbf V}$,
of $\KV$ obviously equals
$\k_v$ since $\A$ has class number $1$.
As we are assuming $\deg v=1$ and $r=2$, we see that $S_v=\Zp$.
As above we find $L_v(\psi,x_v,s_v)=L_v(x_v,s_v)^2$ and
$$L_v(x_v,s_v)=\zeta_{A^\prime,v^\prime}(x_v,2s_v-1)\,.$$
Proposition \ref{wans}, and the results of \cite{w1}, \cite{dv1} now
tell us that the zeroes of $L_v(x_v,s_v)$ are simple, uniquely
determined by their absolute values,  and in
$$\k^\prime_{v^\prime}=\F_2((\sqrt{T}))=\k_v (\Theta)=\k_v(f^\prime)\,$$
for any prime $f\neq v$. So both Conjecture \ref{ffieldGRH2}
and Conjecture \ref{charpsim} are true for $L_v(x_v,s_v)$.
Moreover, if $\lambda$ is one
such zero, one may easily compute $v_{g^\prime}(\lambda)$ which is seen
to be odd. Since the elements of $\F_2((T))$ are precisely the squares
in $\F_2((\sqrt{T}))$ (and so have even valuation), we deduce 
immediately that $\lambda\not\in \k_v=\k_{v,\mathbf V}=\F_2((T))$. 

The same calculation performed at $\infty$ will give an even valuation
and so fails to show that infinitely many zeroes are not
in $\K=\KV$ (though this is indeed likely). 
\end{example}


\begin{thebibliography}{3334}

\bibitem[Am1]{am1} {\sc Y.\ Amice:} Interpolation $p$-adiques,
{\it Bull.\ Soc.\ Math.\ France} {\bf 92} (1964), 117-180.

\bibitem[An1]{a1}  {\sc G.\  Anderson:}  $t$-motives, 
{\it Duke Math. J.} {\bf 53} (1986), 457-502.

\bibitem[BP1]{bp1} {\sc G.\ Boeckle, R.\ Pink:} A cohomological 
theory of crystals over function fields, (in preparation).

\bibitem[Ca1]{ca1} {\sc M.\ Car:} P\`olya's theorem for $\mathbf{F}_q[T]$,
{\it J.\ Number Theory} {\bf 66} (1997), 148-171.

\bibitem[C1]{c1} {\sc L.\ Carlitz:} On certain functions connected with 
polynomials in a Galois field, {\it Duke Math. J.} {\bf 1} (1935), 137--168.

\bibitem[Dr1] {dr1} {\sc V.G.\ Drinfeld:}  
Elliptic modules, Math. Sbornik {\bf 94} 
(1974), 594-627, English transl.:  {\it Math. U.S.S.R. Sbornik} {\bf 23} 
(1976), 561-592.

\bibitem[DV1]{dv1} {\sc J.\ Diaz-Vargas:}  Riemann hypothesis for $\Fp[T]$, 
{\it J.\ Number Theory} {\bf 59} (1996), 313-318.

\bibitem[Go1]{go1} {\sc D.\ Goss}:
{\it Basic Structures of Function Field Arithmetic},
Springer-Verlag, Berlin, 1996.

\bibitem[Go2]{go2} {\sc D.\ Goss}: Some integrals attached to modular forms in
the theory of function fields, in: {\it The Arithmetic of Function Fields} 
(eds:  D.\ Goss et al) de Gruyter (1992), 227-251.

\bibitem[H1] {h1} {\sc D.\ Hayes:} A brief introduction to Drinfeld modules,
in:  {\it The Arithmetic of Function Fields} (eds. D. Goss et al) de Gruyter 
(1992), 1-32.

\bibitem[KS1]{ks1} {\sc N.M.\ Katz, P. Sarnak}:
{\it Random Matrices, Frobenius Eigenvalues and
Monodromy}, Amer. Math. Soc. 1999.

\bibitem[P1]{p1} {\sc R.\ Pink:} The Mumford-Tate conjecture for
Drinfeld-modules, {\it Publ.\ Res.\ Inst.\ Math.\ Sci.\ Kyoto Univ.},
{\bf 33}, (1997), 393-425.

\bibitem[Sh1]{sh1} {\sc J.\ Sheats:} The Riemann hypothesis for the Goss 
zeta function for $\Fq[T]$, {\it J.\ Number Theory} {\bf 71}, (1998),
121-157.

\bibitem[TW1]{tw1} {\sc Y.\ Taguchi, D.\ Wan:} 
$L$-functions of $\varphi$-sheaves and
Drinfeld modules, {\it J. Amer. Math. Soc.} {\bf 9} (1996), 755-781.

\bibitem[TW2]{tw2} {\sc Y.\ Taguchi, D.\ Wan:} Entireness of 
$L$-functions of $\varphi$-sheaves on affine complete 
intersections, {\it J.\ Number Theory} {\bf 63} (1997), 170-179.

\bibitem[W1] {w1} {\sc D.\ Wan:}  On the Riemann hypothesis for the 
characteristic 
$p$ zeta function, {\it J.\ Number Theory} {\bf 58} (1996), 196-212.

\bibitem[Ya1] {ya1} {\sc Z.\ Yang:} Locally analytic functions over
completions of $\Fr[U]$, {\it J.\ Number Theory} {\bf 73}, (1998),
451-458
\end{thebibliography}
\end{document}